\documentclass[10pt]{amsart}
\usepackage{amsfonts,amssymb,amscd}
\usepackage{color}
\input{xy}
\xyoption{all}

\newcommand{\LP}{\operatornamewithlimits{\overrightarrow{\textstyle\prod}}}

\newcommand{\DD}{\mathcal D}
\newcommand{\EE}{\mathcal E}
\newcommand{\F}{\mathcal F}
\newcommand{\LL}{\mathcal L}
\newcommand{\N}{\mathcal N}
\newcommand{\U}{\mathcal U}
\newcommand{\V}{\mathcal V}
\newcommand{\W}{\mathcal W}

\newcommand{\IR}{\mathbb R}

\newcommand{\w}{\omega}
\newcommand{\cbox}{\boxdot}

\newcommand{\id}{\mathrm{id}}

\newcommand{\cl}{\mathrm{cl}}

\newcommand{\ulim}{\mathrm u\mbox{-}\kern-2pt\varinjlim}
\newcommand{\tlim}{\mathrm t\mbox{-}\kern-2pt\varinjlim}
\newcommand{\lclim}{\mathrm{lc}\mbox{-}\kern-2pt\varinjlim}
\newcommand{\llim}{\mathrm{l}\mbox{-}\kern-2pt\varinjlim}

\newcommand{\Int}{\mathrm{Int}}
\newcommand{\tint}{\mathrm{int}}

\newcommand{\glim}{\mathrm{g}\mbox{-}\kern-2pt\varinjlim}

\newtheorem{theorem}{Theorem}
\newtheorem{lemma}{Lemma}

\newtheorem{corollary}{Corollary}
\newtheorem{proposition}{Proposition}

\theoremstyle{definition}

\newtheorem{remark}{Remark}
\newtheorem*{assumption}{Assumption}

\title[]{Diffeomorphism groups of non-compact manifolds \\
endowed with the Whitney $C^\infty$-topology}

\author[T. Banakh]{Taras Banakh}
\address[T. Banakh]{Ivan Franko National University of Lviv (Ukraine) and Jan Kochanowski University in Kielce (Poland)}
\email{t.o.banakh@gmail.com}

\author[T. Yagasaki]{Tatsuhiko Yagasaki}
\address[T. Yagasaki]{Graduate School of Science and Technology, Kyoto Institute of Technology, Kyoto, 606-8585, Japan}
\email{yagasaki@kit.ac.jp}
\thanks{The first author has been partially financed by NCN means granted by decision DEC-2011/01/B/ST1/01439. The second author was supported by JSPS Grant-in-Aid for Scientific Research (No.22540081).}

\subjclass[2010]{57N20; 57S05; 58D05; 46A13}
\keywords{Diffeomorphism group, the Whitney topology, $\sigma$-compact manifold, LF-space}

\begin{document}
\baselineskip 5mm

\begin{abstract} Suppose $M$ is a non-compact connected $n$-manifold without boundary, $\DD(M)$ is the group of $C^\infty$-diffeomorphisms of $M$
endowed with the Whitney $C^\infty$-topology and $\DD_0(M)$ is the identity connected component of $\DD(M)$, which is an open subgroup in the group $\DD_c(M)\subset\DD(M)$ of compactly supported diffeomorphisms of $M$.
It is shown that $\DD_0(M)$ is homeomorphic to $N \times \IR^\infty$ for an $l_2$-manifold $N$ whose topological type is uniquely determined by the homotopy type of $\DD_0(M)$.
For instance, $\DD_0(M)$ is homeomorphic to $l_2 \times \IR^\infty$
if $n = 1, 2$ or $n = 3$ and $M$ is orientable and irreducible.
We also show that for any compact connected $n$-manifold $N$ with non-empty boundary $\partial N$ the group $\DD_0(N\setminus\partial N)$ is homeomorphic to $\DD_0(N; \partial N) \times \IR^\infty$, where $\DD_0(N;\partial N)$ is the identity component of the group $\DD(N;\partial N)$ of diffeomorphisms of $N$ that do not move points of the boundary $\partial N$.
\end{abstract}

\maketitle

\section{Introduction}

In this paper we continue the study of the topological structure of diffeomorphism groups of non-compact smooth manifolds endowed with the Whitney $C^\infty$-topology, started in \cite{BMSY1}.
Suppose $M$ is a $\sigma$-compact smooth $n$-manifold without boundary.
Let $\DD(M)$ denote the group of diffeomorphisms of $M$ endowed with the Whitney $C^\infty$-topology
(called the very-strong $C^\infty$-topology in \cite{Illman})
and $\DD_0(M)$ be the identity connected component of $\DD(M)$.
The group $\DD(M)$ includes the normal subgroup $\DD_c(M)$
consisting of diffeomorphisms with compact support.

In \cite[Theorem 4, Theorem 6.8]{BMSY1} we have shown that $\DD_c(M)$ is a paracompact $(l_2 \times \IR^\infty)$-manifold and $\DD_0(M)$ is an open subgroup of $\DD_c(M)$, while $\DD(M)$ itself is locally homeomorphic to the box product $\square_{i \in \w} l_2$.
Here $l_2$ is the separable Hilbert space and $\IR^\infty$ is the direct limit of the sequence $(\IR^n)_{n \in \w}$, where
$\IR^n$ is identified with the hyperspace $\IR^n \times \{ 0 \}$ in $\IR^{n+1}$.

In a series of papers \cite{BR1, BR2, BR3} T.~Banakh and D.~Repov\v{s} studied topological properties of direct limits in the category of uniform spaces.
These results were applied in \cite{BMRSY} to yield a simple criterion for recognizing topological groups homeomorphic to open subspaces of $l_2 \times \IR^\infty$ (see Theorem~\ref{thm_BR} in Section 2).
In this paper we apply this criterion to obtain the following important conclusion on the group $\DD_c(M)$.

\begin{theorem}\label{thm_diff} For any non-compact $\sigma$-compact smooth $n$-manifold $M$ without boundary the group $\DD_c(M)$ is homeomorphic to  an open subspace of $l_2 \times \IR^\infty$.
\end{theorem}

In \cite{MS} K. Mine and K. Sakai proved Triangulation Theorem for open subsets of $l_2 \times \IR^\infty$ (see Theorem~\ref{thm_MS} in Section 2): any open subset $U$ of $l_2 \times \IR^\infty$ is homeomorphic to $N \times {\Bbb R}^\infty$ for some $l_2$-manifold $N$ whose topological type is uniquely determined by the homotopy type of $U$. Combining this result with Theorem~\ref{thm_diff}, we obtain:

\begin{corollary}\label{cor_factor}
The group $\DD_0(M)$ is homeomorphic to $N \times {\Bbb R}^\infty$ for some $l_2$-manifold $N$ whose topological type is uniquely determined by the homotopy type of $\DD_0(M)$. In particular, if $\DD_0(M)$ is homotopy equivalent to an $l_2$-manifold $N$, then $\DD_0(M)$ is homeomorphic to $N \times {\Bbb R}^\infty$.
\end{corollary}

In some specific cases we can detect the homotopy and topological types of $\DD_0(M)$. Below the symbol $\approx$ denotes the topological equivalence.

\begin{corollary}\label{cor_top_type}
Let $M$ be a non-compact connected smooth $n$-manifold without boundary.
\begin{enumerate}
\item
If $1\le n\le 2$, then $\DD_0(M) \approx l_2 \times \IR^\infty$.
\item
If $n = 3$ and the manifold $M$ is orientable and irreducible
 (i.e., any smooth $2$-sphere in $M$ bounds a $3$-ball in $M$),
 then $\DD_0(M) \approx l_2 \times \IR^\infty$.
\item
If $M$ is diffeomorphic to the interior $N \setminus \partial N$
 of a compact connected smooth $n$-manifold $N$ with boundary, then
$\DD_0(M) \approx \DD_0(N, \partial N) \times \IR^\infty$.
In particular, $\DD_0(M) \approx l_2 \times \IR^\infty$ if $\DD_0(N;\partial N)$ is contractible.
\end{enumerate}
\end{corollary}

In this corollary $\DD_0(N;\partial N)$ stands for the identity connected component of the group $\DD(N;\partial N)$ of diffeomorphisms of $N$ that do not move points of the boundary $\partial N$. Because of the compactness of the manifold $N$, the Whitney $C^\infty$-topology on $\DD(N;\partial N)$ is metrizable and coincides with the compact-open $C^\infty$-topology.
Corollary~\ref{cor_top_type} implies, for instance, that $\DD_0(M) \approx l_2 \times \IR^\infty$ if $M$ is the 3-dimensional Euclidean space ${\Bbb R}^3$
 or the Whitehead contractible 3-manifold \cite{Hempel}.

The authors hope that the results in \cite{BMSY1} and in this paper will clarify ambiguity in literatures on the Whitney $C^\infty$-topology on diffeomorphism groups of non-compact manifolds.

\section{Open subspaces of LF-spaces}

The study of the topological structure of topological groups has a long history (cf. \cite{BMRSY}).
Topological groups locally homeomorphic to the separable Hilbert space $l_2$ were characterized by Dobrowolski and Toru\'nczyk \cite{DT};
A topological group $G$ is a separable (finite or infinite dimensional) Hilbert manifold  if and only if
$G$ is a Polish ANR.
Here,
a Polish space means a separable completely metrizable space and an ANR means
an absolute neighborhood retract for metric spaces.
In this paper every (finite or infinite dimensional) manifold is assumed to be paracompact.

In this preliminary section we recall a criterion for recognizing topological groups homeomorphic to open subsets of the separable non-metrizable LF-space $l_2 \times \IR^\infty$.  First we recall some necessary definitions.

Suppose $G$ is a topological group with the neutral element $e$.
A tower of subgroups of $G$ means a sequence $(G_n)_{n \in \w}$  of subgroups of $G$ such that
$$G_0\subset G_1\subset G_2\subset\cdots \hspace{5mm} and \hspace{5mm} G=\bigcup_{n\in\w}G_n.$$
Following \cite{BMRSY}, we say that $G$ carries the {\it strong topology} with respect to the tower $(G_n)_{n \in \w}$
if for any neighborhood $U_n$ of the neutral element $e$ in $G_n$, $n \in \w$,  the group product
$$\LP_{n\in\w}U_n=\bigcup_{n\in\w} U_0U_1\cdots U_n$$is a neighborhood of $e$ in $G$. In this case the topology of $G$ coincides with the topology of the direct limit $\glim G_n$ of the tower $(G_n)_{n\in\w}$ in the category of topological groups, which means that $G$ carries the strongest group topology such that the inclusion maps $G_n\to G$, $n\in\w$, are continuous.

A subgroup $H$ of a topological group $G$ is called {\em locally topologically complemented} (LTC) in $G$ if $H$ is closed in $G$ and the quotient map $q : G \to G/H=\{xH:x\in G\}$ is a locally trivial bundle. This condition is equivalent to saying that $q$ has a local section at some point of $G/H$.
Here, a {\em local section} of a map $q: X \to Y$ at a point $y \in Y$ means a continuous map $s: U \to X$ defined on a neighborhood $U$ of $y$ in $Y$ such that $q\circ s=\id_U$.

A closed subset $A$ of a topological space $X$ is called a ({\em strong}) {\em $Z$-set in} $X$ if for any open cover $\U$ of $X$ there is a continuous map $f:X\to X$ such that $f$ is $\U$-near to the identity $\id_X:X\to X$ and (the closure of) the set $f(X)$ does not intersect $A$.
A point $x$ of $X$ is called a ({\em strong}) {\em $Z$-point} if the singleton $\{x\}$ is a (strong) $Z$-set in $X$.
Every point in an infinite-dimensional Hilbert manifold is a strong $Z$-point.

The following criterion is obtained in \cite[Theorem 1.6]{BMRSY}.

\begin{theorem}[Banakh-Mine-Repov\v{s}-Sakai-Yagasaki]\label{thm_BR}
A non-metrizable topological group $G$ is homeomorphic to an open subset of\, $\IR^\infty$ or $l_2 \times \IR^\infty$ if $G$ carries the strong topology with respect to a tower $(G_n)_{n\in\w}$ of closed subgroups of $G$
such that for each $n \in \w$
{\rm (i)} $G_n$ is a separable Hilbert manifold, {\rm (ii)} $G_n$ is LTC in $G_{n+1}$ and
{\rm (iii)} each $Z$-point of the quotient space $G_{n+1}/G_n$ is a strong $Z$-point.
\end{theorem}

Open subspaces of $l_2\times\IR^\infty$ were studied in \cite{MS} and the following Triangulation Theorem was obtained.

\begin{theorem}[Mine-Sakai]\label{thm_MS}
{\rm (1)} Each open subspace $X$ of $l_2\times\IR^\infty$ is homeomorphic to the product $K\times l_2\times\IR^\infty$
for a locally finite simplicial complex $K$.

{\rm (2)} Two open subspaces of $l_2\times\IR^\infty$ are homeomorphic if and only if they are homotopically equivalent.
\end{theorem}

Note that the product $N = K \times l_2$ is an $l_2$-manifold and its topological type is determined by its homotopy type.

In the next section we apply Theorem~\ref{thm_BR} to the diffeomorphism groups of non-compact manifolds.
To check the conditions imposed in this criterion we need some results on small box products.
For a sequence of pointed spaces $(X_i, \ast_i)$, $i \in \w$,
the box product $\square_{i \in \w} X_i$ is the Cartesian product of $X_i$, $i \in \w$,  endowed with the box topology. This topology is generated
 by the base consisting of boxes $\square_{i\in\w}U_i$, where $U_i$ is an open set of $X_i$.
The small box product $\cbox_{i \in \w} X_i$ is the subspace of $\square_{i \in \w} X_i$
defined by
$$\cbox_{i\in\w} X_i = \big\{(x_i)_{i\in\w}\in \square_{i\in\w} X_i :
 \exists \,k\in\w\;\forall\, i\ge k, \ x_i= \ast_i \,\big\}.$$
The subspace topology on $\cbox_{i\in\w} X_i$ is generated
 by the base consisting of small boxes $\cbox_{i\in\w}U_i$, where $U_i$ is an open set of $X_i$.
It is known that $\cbox_{i\in\w} \IR \approx \IR^\infty$ and $\cbox_{i\in\w} l_2 \approx l_2 \times \IR^\infty$, and in \cite[Theorem 4]{BMSY1} we have shown that
for a non-compact $\sigma$-compact smooth $n$-manifold $M$ without boundary, the pair $(\DD(M), \DD_c(M))$ is locally homeomorphic to the pair $(\square_{i\in\w} l_2, \cbox_{i\in\w} l_2)$ at $\id_M$.

Suppose $G$ is a topological group with the neutral element $e$ and $(G_n)_{n \in \w}$ is a tower of subgroups of $G$.
This tower induces the small box product $\cbox_{i \in \w} G_i$ and
the left multiplication map
$$p : \cbox_{i \in \w} G_i \to G, \ p(x_0,\dots,x_k,e,e,\dots) = x_0 \cdot x_1\cdots x_k$$
The map $p$ is continuous (\cite[Lemma 2.10]{BMSY1}), and
the group $G$ carries the strong topology with respect to the tower $(G_n)_{n \in \w}$ if and only if
the map $p$ is open at $(e)_{i \in \w}  \in \cbox_{i \in \w} G_i$.

\begin{remark}\label{rmk_op}
If the map $p$ has a local section $s$ at the neutral element $e \in G$, then
\begin{itemize}
\item[(i)\,] for any $x \in G$ and any $\boldsymbol{x} \in p^{-1}(x)$
one can modify the local section $s$ to obtain
a local section $\sigma$ of $p$ at $x$ with $\sigma(x) = \boldsymbol{x}$ and hence
\item[(ii)] the map $p$ is open, which implies that $G$ carries the strong topology with respect to the tower $(G_n)_{n \in \w}$.
\end{itemize}
\end{remark}

We close this section with some remarks on LTC subgroups.

\begin{lemma}\label{lem_ltc}
Suppose $G$ is a topological group and $K \subset H$ are closed subgroups of $G$.
\begin{itemize}
\item[(1)] If $G$ is metrizable, then so is the quotient space $G/H$.
\item[(2)] If $K$ is LTC in $H$ and $H$ is LTC in $G$, then $K$ is LTC in $G$.
\item[(3)] If $H$ is LTC in $G$, then the map $\pi :  G/K \to G/H$, $\pi(gK) = gH$,
 is a locally trivial bundle with the fiber $H/K$.
\end{itemize}
\end{lemma}

\begin{proof}
(1) The group $G$, being metrizable, admits a right invariant metric $d$ generating the topology of $G$.
Then the topology of the quotient space $G/H$ is generated by
the metric $\rho$ defined by
$$\rho(xH,yH)=\inf\{d(a,b):a\in xH,\; b\in yH\} \hspace{5mm} \mbox{for \ }xH,yH\in G/H.$$

(2)
By the assumption,
the projection $G \to G/H$ has a local section
 $\sigma : (W, \overline{e}) \to (G, e)$ at $\overline{e} \in G/H$ and
the projection $H \to H/K$ also has a local section
 $\tau : U \to H$ at $\overline{e} \in H/K$.
 Consider the projection $\pi : G/K \to G/H$.
The map $\sigma$ determines the map
 $\sigma_0 : \pi^{-1}(W) \to H/K$, $\sigma_0(x) = \sigma(\pi(x))^{-1}x$.
Since $\sigma_0(\overline{e}) = \overline{e} \in U$,
there exists an open neighborhood $V$ of $\overline{e}$ in $\pi^{-1}(W)$
 such that $\sigma_0(V) \subset U$.
The required local section $s : V \to G$ for the projection $G \to G/K$ is defined by
 $s(x) = \sigma(\pi(x)) \tau(\sigma_0(x))$.

(3)
The projection $G \to G/H$ has a local section
 $\sigma : W \to G$ at any point $x_0 \in G/H$.
The associated trivialization
 $\phi : W \times H/K \approx \pi^{-1}(W)$ is defined by
 $\phi(x, y) = \sigma(x)y$.
\end{proof}

%%%%%%%%%%%%%%%

\section{Spaces of embeddings and Bundle Theorem}\label{Diffeo groups}

Suppose $M$ is a smooth $n$-manifold without boundary.
The diffeomorphism group $\DD(M)$ is endowed with the Whitney $C^\infty$-topology and
any subset $K$ of $M$ induces a closed subgroup $\DD(M; K) = \{ h \in \DD(M) : h|_{K} = \id_{K}\}$.
Let $\DD_0(M, K)$ denote the identity connected components of $\DD(M, K)$.

For a smooth submanifold $L$ of $M$ and a subset $K$ of $L$
let $\EE_K(L, M)$ denote the space of $C^\infty$-embeddings $f : L \to M$ with $f|_K = \id_K$.
This space and its subspaces are endowed with the compact-open $C^\infty$-topology.
There is a natural restriction map
$$r : \DD(M;K) \to \EE_K(L, M), \quad r(h) = h|_L.$$
Let $\EE_K^{\star}(L, M)$ denote the  image of $r$, i.e., ${\rm Im}\,r = \{ h|_L : h \in \DD(M, K)\}$.
This is the subspace of $\EE_K(L, M)$ consisting of extendable $C^\infty$-embeddings.
Note that the group $\DD(M, K)$ acts continuously on $\EE_K(L, M)$ by the left composition and under this action
the map $r$ is the orbit map at the inclusion $i_L : L \subset M$ and $\EE_K^{\star}(L, M)$ is the orbit of $i_L$.
Let $\EE_K(L, M)_0$ denote the connected component of the inclusion $i_L$ in $\EE_K(L, M)$.
Then the map $r$ induces the restriction map
$$r_0 : \DD_0(M;K) \to \EE_K(L, M)_0, \quad r_0(h) = h|_L.$$

The followings are the classical bundle theorem in codimension 0 and its complements (cf.\ \cite{Cerf},  \cite{Hmlt} \cite{Palais1}, \cite{Seeley}).
Below we impose the following condition on $M$ and $K \subset L$.

\begin{assumption}
We assume that $M$ is a $\sigma$-compact smooth $n$-manifold without boundary and
$K \subset L$ are smooth $n$-submanifolds of $M$ such that
$K, L$ are closed subsets of $M$, $K \subset \Int L$ and $\cl_M(L \setminus K)$ is compact.
\end{assumption}

\begin{theorem}\label{thm_bdl} \mbox{}
For any closed subset $C$ of $M$ with $C \cap L = \emptyset$ there exists a neighborhood $\U$ of
the inclusion $i_L$ in $\EE_K(L, M)$ and a map
$s : \U \to \DD_0(M;K \cup C) \subset \DD(M, K)$ such that $s(i_L) = \id_M$ and
$s(g)|_L = g$ for $g \in \U$.
\end{theorem}

The map $s$ is exactly a local section of the restriction map $r : \DD(M, K) \to \EE_K(L, M)$ at $i_L$.

\begin{remark}\label{rmk_co-w} In some literatures, the group $\DD_0(M;K \cup C)$ is endowed with the compact-open $C^\infty$-topology.
However, since $\cl_M(L \setminus K)$ is compact, we can enlarge $C$ so that
$\cl_M(M - (K \cup C))$ is compact.
In this case the Whitney $C^\infty$-topology on the group $\DD_0(M;K \cup C)$ coincides with the compact-open $C^\infty$-topology.
\end{remark}

Theorem~\ref{thm_bdl} extends to the following form.

\begin{corollary}\label{cor_bdl-1}
For any $f \in \EE_K(L, M)$ there exist an open neighborhood $\V_f$ of $f$ in $\EE_K(L, M)$ and a map $$\eta_f : \V_f \to \DD_0(M, K)$$
such that $\eta_f(f) = \id_M$ and $\eta_f(g)f = g$ for $g \in \V_f$. Furthermore, the following holds:
\begin{itemize}
\item[(i)\,] If $\V_f \cap \EE_K^\star(L, M) \neq \emptyset$, then $\V_f \subset \EE_K^{\star}(L, M)$
and the map $r$ has a section on $\V_f$. Hence, for any $f \in \EE_K^\star(L, M)$ the map $r$ has a local section at $f$.
\item[(ii)] If $\V_f \cap {\rm Im}\,r_0 \neq \emptyset$, then
$\V_f \subset {\rm Im}\,r_0$ and the map $r_0$ has a section on $\V_f$.
Hence, for any $f \in {\rm Im}\,r_0$ the map $r_0$ has a local section at $f$.
\end{itemize}
\end{corollary}

\begin{corollary}\label{cor_bdl-2} {\rm (1)} {\rm (i)} $\EE^{\star}_K(L, M)$ is a clopen subset of $\EE_K(L, M)$.
{\rm (ii)} ${\rm Im}\,r_0 = \EE_K(L, M)_0$.
\begin{enumerate}
\item[(2)]
\begin{itemize}
\item[(i)\,] The map $r : \DD(M, K) \to \EE_K^\star(L, M)$ is a topological principal bundle with structure group $\DD(M, L)$.
\item[(ii)] The map $r_0 : \DD_0(M, K) \to \EE_K(L, M)_0$ is a topological principal bundle with structure group $\DD(M, L) \cap \DD_0(M, K)$.
\end{itemize}
\item[(3)]
\begin{itemize}
\item[(i)\,] $\DD(M, L)$ is LTC in $\DD(M, K)$ and $\DD(M, K)/\DD(M, L) \approx \EE_K^\star(L, M)$.
\item[(ii)] $\DD(M, L) \cap \DD_0(M, K)$ is LTC in $\DD_0(M, K)$ and
$\DD_0(M, K)\big/\big(\DD(M, L) \cap \DD_0(M, K)\big) \approx \EE_K(L, M)_0$.
\end{itemize}
\end{enumerate}
\end{corollary}

The next theorem is also a classical one (cf.\ \cite{Les}).

\begin{theorem}\label{thm_Frechet}
If $L \setminus K \neq \emptyset$, then the spaces
$\DD(M, K \cup (M \setminus L))$ and $\EE_K(L, M)$ are infinite-dimensional separable Fr\'echet manifolds and hence they are topological $l_2$-manifolds.
\end{theorem}

Since $\DD_0(M, K \cup (M \setminus L))$ and $\EE^{\star}_K(L, M)$, $\EE_K(L, M)_0$ are open subspaces of
$\DD(M, K \cup (M \setminus L))$ and $\EE_K(L, M)$ respectively, they are
also infinite-dimensional separable Fr\'echet manifolds and hence are topological $l_2$-manifolds.

For the sake of completeness we include a sketch of the proofs of Corollaries~\ref{cor_bdl-1} and~\ref{cor_bdl-2}.

\begin{proof}[\bf Proof of Corollary~\ref{cor_bdl-1}]
We apply Theorem~\ref{thm_bdl} to the pair $K \subset f(L)$ so to obtain
an open neighborhood $\V$ of $i_{f(L)}$ in $\EE_K(f(L), M)$ and
a map $s' : \V \to \DD_0(M, K)$ such that $s'(i_{f(L)}) = \id_M$ and $s'(g')|_{f(L)} = g'$ for $g' \in \V$.
The diffeomorphism $f : L \cong f(L)$ induces the homeomorphism $f^\ast : \EE_K(f(L), M) \approx \EE_K(L, M)$.
Then $\V_f = f^\ast(\V)$ is a neighborhood of $f^\ast(i_{f(L)}) = f$ in $\EE_K(L, M)$ and
the required map $\eta_f : \V_f \to \DD_0(M, K)$ is defined by
$\eta_f(g) = s'(gf^{-1})$ for $g \in \V_f$.

The additional statements are verified as follows:
(i) If $\V_f \cap \EE_K^\star(L, M) \neq \emptyset$, then
we can find $f' \in \V_f \cap \EE_K^\star(L, M)$ and $h \in \DD(M, K)$ with $h|_L = f'$.
The desired section $s_f : \V_f \to \DD(M, K)$ of $r$
is defined by $s_f(g) = \eta_f(g) (\eta_f(f'))^{-1} h$ for $g \in \V_f$.
(ii) If $\V_f \cap {\rm Im}\,r_0 \neq \emptyset$, then we can find
$f' \in \V_f \cap {\rm Im}\,r_0$ and $h \in \DD_0(M, K)$ with $h|_L = f'$.
The desired section $s_f : \V_f \to \DD_0(M, K)$ of $r_0$ is defined by $s_f(g) = \eta_f(g) (\eta_f(f'))^{-1} h$ for $g \in \V_f$.
\end{proof}

\begin{proof}[\bf Proof of Corollary~\ref{cor_bdl-2}]
(1) (i) For any $f \in cl\,\EE^{\star}_K(L, M)$ we can find $\V_f$ and $\eta_f$ as in Theorem~\ref{cor_bdl-1}.
Since $\V_f \cap \EE^{\star}_K(L, M) \neq \emptyset$, by Theorem~\ref{cor_bdl-1} (1) we have
$f \in \V_f \subset \EE^{\star}_K(L, M)$.
This means that $\EE^{\star}_K(L, M)$ is clopen.

(ii) For any $f \in cl\,{\rm Im}\,r_0$ we can find $\V_f$ and $\eta_f$ as in Theorem~\ref{cor_bdl-1}.
Since $\V_f \cap {\rm Im}\,r_0 \neq \emptyset$, by Theorem~\ref{cor_bdl-1} (2) we have
$f \in \V_f \subset {\rm Im}\,r_0$.
This means that ${\rm Im}\,r_0$ is clopen.
Since $\DD_0(M, K)$ is connected, so is ${\rm Im}\,r_0$, which implies the conclusion.

(2) (i)
The group $\DD(M, L)$ acts continuously on $\DD(M, K)$ by the right composition and
this action preserves the fibers of the map $r$.
By Corollary~\ref{cor_bdl-1}\,(i) for any $f_0 \in \EE_K^\star(L, M)$ the map $r : \DD(M, K) \to \EE_K^\star(L, M)$ admits a local section $s : \U \to \DD(M, K)$ at $f_0$.
Then a $\DD(M, L)$-equivariant trivialization $\phi : r^{-1}(\U) \to  \U \times \DD(M, L)$ of the map $r : r^{-1}(\U) \to \U$ is obtained by $\phi(g) = (r(g), s(r(g))^{-1}g)$.
The inverse homeomorphism $\psi = \phi^{-1} : \U \times \DD(M, L) \to r^{-1}(\U)$  is given by $\psi(f, h) = s(f)h$.

(ii) The group $\DD(M, L) \cap \DD_0(M, K)$ acts continuously on $\DD_0(M, K)$ by the right composition and
this action preserves the fibers of the map $r_0$.
By (1)(ii) and Corollary~\ref{cor_bdl-1}\,(ii) for any $f_0 \in \EE_K(L, M)_0$ the map $r_0 : \DD_0(M, K) \to \EE_K(L, M)_0$ admits a local section $s : \U \to \DD_0(M, K)$ at $f_0$.
Then a $\big(\DD(M, L) \cap \DD_0(M, K)\big)$-equivariant trivialization $\phi : r_0^{-1}(\U) \to \U \times \big(\DD(M, L) \cap \DD_0(M, K)\big)$ of the map $r_0 : r_0^{-1}(\U) \to \U$ is obtained by $\phi(g) = (r_0(g), s(r_0(g))^{-1}g)$.
The inverse homeomorphism $\psi = \phi^{-1} : \U \times \big(\DD(M, L) \cap \DD_0(M, K)\big) \to r_0^{-1}(\U)$ is given by $\psi(f, h) = s(f)h$.

(3) By (2)(i) the map $r$ has the factorization
$$\xymatrix@M+2pt{
& \DD(M, K) \ar[dl]_{\mbox{$\pi$}} \ar[dr]^{\mbox{$r$}} & & \\
\DD(M, K)/\DD(M, L) \ar[rr]^{\ \ \ \ \ \mbox{$\phi$}\ }_{\ \ \ \ \approx} & & \EE_K^{\star}(L, M),
}$$
where the homeomorphism $\phi$ is defined by $\phi(h\,\DD(M, L)) = h|_L$.
This implies the assertion (i).
Similarly, from (2)(ii) follows the statement (ii).
\end{proof}

%%%%%%%%%%%%%%%

\section{Diffeomorphism groups of non-compact manifolds}\label{Diffeo groups}

Suppose $M$ is a non-compact $\sigma$-compact smooth $n$-manifold without boundary.
We can represent $M$ as the countable union $M=\bigcup_{i\in\w}M_i$ of compact $n$-submanifolds $M_i$, $i\in\w$, of $M$ such that $M_0 \neq \emptyset$ and $M_i \subsetneqq \Int M_{i+1}$ for all $i\in \w$.
Let $M_{-1}=\emptyset$ and consider the $n$-submanifolds $K_i = M \setminus \Int M_i$, $i\in\w$, of $M$.
Then we obtain the group $G = \DD_c(M)$ and
the tower $G_i = \DD(M; K_i)$, $i\in\w$, of closed subgroups of $G$.
This tower induces the small box product $\cbox_{i \in \w} G_i$ and
the left multiplication map $p : \cbox_{i \in \w} G_i \to G$.

We shall show that the tower $(G_i)_{i \in \w}$ has the properties listed in Proposition~\ref{prop_G_i} below.
Hence, Theorem~\ref{thm_diff} now follows from Theorem~\ref{thm_BR}.

\begin{proposition}\label{prop_G_i}
{\rm (1)}
The group $G$ is paracompact, but not metrizable.
\begin{itemize}
\item[{\rm (2)}] For each $i \in \w$
\begin{itemize}
\item[{\rm (i)}\ ] $G_i$ is a separable $l_2$-manifold,
\item[{\rm (ii)}\,] $G_i$ is LTC in $G_{i+1}$ and
\item[{\rm (iii)}] $G_{i+1}/G_i$ is a separable $l_2$-manifold.
\end{itemize}
\item[{\rm (3)}] The left multiplication map $p : \cbox_{i \in \w} G_i \to G$
admits a local section at the identity $\id_M$.
Hence, the map $p$ is an open map.

\item[{\rm (4)}] The group $G$ carries the strong topology with respect to the tower $(G_i)_{i \in \w}$.
Hence $G = \glim G_i$.
\end{itemize}
\end{proposition}

Let $H$ and $H_i$, $i \in \w$, denote the identity connected components of $G$ and $G_i$, $i \in \w$, respectively.

\begin{proposition}\label{prop_H_i}  
{\rm (1)$'$} \, $H$ is an open subgroup of $G$ and the sequence $(H_i)_{i\in\w}$ forms a tower of of closed subgroups 
\begin{itemize}
\item[] of $H$.
\item[{\rm (2)$'$}] The tower $(H_i)_{i \in \w}$ has the properties $(1) - (4)$ listed in Proposition~\ref{prop_G_i} with respect to the group $H$.
\end{itemize}
\end{proposition}

\begin{proof}[\bf Proof of Proposition~\ref{prop_G_i}]
(1) The paracompactness of $G$ is shown in (the proof of) \cite[Theorem 6.8]{BMSY1} (cf. \cite[Proposition 4.1(2))]{BMSY1}.
Using the non-compactness of $M$ and the diagonal argument, we can easily show that $G$ is not first countable.
This implies the non-metrizability of $G$.

(2) (i) Since $M\setminus\Int\,K_i = M_i$ is compact, the group
$G_i$ is an infinite-dimensional separable Fr\'echet manifold (Theorem~\ref{thm_Frechet}).

(ii), (iii) We apply Corollary~\ref{cor_bdl-2}\,(3)(i) to $K_{i+1} \subset K_i$ to conclude that
$G_i$ is LTC in $G_{i+1}$ and $G_{i+1}/G_i \approx \EE^{\star}_{K_{i+1}}(K_i, M)$.
The latter is an infinite-dimensional separable Fr\'echet manifold (Theorem~\ref{thm_Frechet}).

(3), (4) The existence of a local section of the map $p$ follows from
\cite[Proposition 5.5\,(2)]{BMSY1} (cf. Proof of \cite[Theorem 6.8]{BMSY1}),
which treats a general pair of a transformation group and
the subgroup of those with compact support endowed with a Whitney-like topology.
Since its full-generality seems to prevent from understanding the essence for our purpose here,
we include below a short self-contained proof specified to the case of $G = \DD_c(M)$.
The remaining assertions in the statements (3), (4) follow from Remark~\ref{rmk_op}.

Below we use the following notations:
For a subset $K$ of $M$ let $G_K = \DD_c(M, K)$ and $G(K) = \DD_c(M, M - K)$.
Then $G_i = G_{K_i} = G(M_i)$. Let
$$\mbox{$F_i = M_{2i} - {\rm Int}\,M_{2i-1}$ \quad and \quad
$L_i = M_{2i+1} - {\rm Int}\,M_{2i}$ \quad $(i \in \w)$.}$$
There exists a sequence $(N_i)_{i\in\w}$ of compact $n$-submanifolds of $M$
 such that $L_i \subset \tint_M N_i$, $N_i \subset \tint_M M_{2i+2} \setminus M_{2i-1}$ and
$N_i \cap N_j = \emptyset$.
Then
$$\mbox{$\F = (F_i)_{i\in\w}$, $\LL = (L_i)_{i\in\w}$ \ and \ $\N = (N_i)_{i\in\w}$}$$
are discrete families of compact $n$-submanifolds of $M$.
Let $F = \bigcup_{i\in\w} F_i$, $L = \bigcup_{i\in\w} L_i$ and $N = \bigcup_{i\in\w} N_i$.
We have $M = F \cup L$ and $G_L = G(F)$.

(i) First consider the map $\lambda_\F : \cbox_{i\in\w} G(F_i) \to G(F)$ defined by
$$\lambda_\F(h_0, h_1, \cdots, h_m, \id_M, \cdots) = h_0h_1h_2 \cdots h_m.$$
The map $\lambda_\F$ is an open embedding. In fact, since $\F = (F_i)_{i\in\w}$ is discrete, the map $\lambda_\F$ has the natural extension
$\widetilde{\lambda}_\F : \square_{i\in\w}G(F_i) \to \DD(M; M - F)$ defined by
$$\widetilde{\lambda}_\F((h_i)_{i\in\w})|_{M \setminus F} = \id \hspace{5mm} \text{and} \hspace{5mm} \widetilde{\lambda}_\F((h_i)_{i\in\w})|_{F_j} = h_j|_{F_j} \mbox{ \ for \ }j \in \w
.$$
Due to the definition of the Whitney $C^\infty$-topology, it is seen that $\widetilde{\lambda}_\F$ is an open embedding.
Since $\widetilde{\lambda}_\F^{-1}(G(F)) = \cbox_{i\in\w}G(F_i)$, the map $\lambda_\F$ is also an open embedding.
The map $\lambda_\N: \cbox_{i\in\w} G(N_i) \to G(N)$ is defined similarly.

(ii) Next we show that the map
$$\theta : \cbox_{i\in\w} G(N_i) \times G(F) \longrightarrow G : \ \theta((g_i)_{i\in\w}, h) = \lambda_{\N}((g_i)_{i\in\w}) h$$
has a local section at $\id_M$.
For each $i \in \w$, applying Theorem~\ref{thm_bdl} to $C = M - {\rm Int}\,N_i$, we find
an open neighborhood $\V_i$ of the inclusion $i_{L_i}$ in $\EE^\star(L_i, M)$
and a map $s_i : \V_i \to G(N_i)$ such that $s_i(f)|_{L_i} = f$ for each $f \in \V_i$ and $s_i(i_{L_i}) = \id_M$.
The maps $s_i$, $i \in \w$, determine the map
$$s = \cbox_{i\in\w} s_i : \cbox_{i\in\w} \V_i \longrightarrow \cbox_{i\in\w} G(N_i) : s((f_i)_{i\in\w}) = (s_i(f_i))_{i\in\w}.$$
Consider the map
$$r_\LL : G \longrightarrow \cbox_{i\in\w} \EE^\star(L_i, M) : r_\LL(g) = (g|_{L_i})_{i\in\w}.$$
The continuity of this map also relies on the definition of the Whitney $C^\infty$-topology.
The inverse image $\V = r_\LL^{-1}(\cbox_{i\in\w} \V_i)$ is an open neighborhood of $\id_M$ in $G$ and
we obtain the composition
$$\eta = \lambda_\N s \,r_\LL : \V \to G(N).$$
Note that $\eta(g)^{-1}g \in G_{L} = G(F)$ for each $g \in \V$, since
$$\eta(g) = \lambda_\N s(g|_{L_i})_i = \lambda_\N (s_i(g|_{L_i}))_i \quad
\text{and} \quad
\eta(g)|_{L_i} = s_i(g|_{L_i})|_{L_i} = g|_{L_i}.$$
The desired local section of $\theta$ at $\id_M$ is defined by
$$\sigma_0 : \V \longrightarrow \cbox_{i\in\w} G(N_i) \times G(F) : \ \sigma_0(g) = (sr_\LL(g), \eta(g)^{-1}g).$$
In fact, for any $g \in \V$ it follows that
$$\theta\sigma_0(g) = \theta(sr_\LL(g), \eta(g)^{-1}g)
= \lambda_\N sr_\LL(g) \eta(g)^{-1} g  = \eta(g) \eta(g)^{-1} g   = g.$$

(iii) Consider the map
$$\rho = \theta (\id \times \lambda_\F) : \cbox_{i\in\w} G(N_i) \times \cbox_{i\in\w} G(F_i) \longrightarrow G : \
\rho((g_i)_{i\in\w},(h_i)_{i\in\w})
= \lambda_\N((g_i)_{i\in\w})\lambda_\F((h_i)_{i\in\w}).$$
Since $\lambda_{\F}$ is an open embedding,
the image $\W = \cbox_{i\in\w} G(N_i) \times {\rm Im}\, \lambda_{\F}$ is an open neighborhood of $\sigma_0(\id_M) = ((\id_M)_{i\in\w}, \id_M)$
 in $\cbox_{i\in\w} G(N_i) \times G_{L}$. Hence,
if we replace $\V$ by a smaller one $\U$, then $\sigma_0(\U) \subset \W$ and
a local section of $\rho$ on $\U$ is defined by
$$\sigma = (\id \times \lambda_\F^{-1}) \sigma_0 : \U \longrightarrow \cbox_{i\in\w} G(N_i) \times \cbox_{i\in\w} G(F_{i}).$$

For each $h \in \U$
the image $\sigma(h) = ((f_i)_{i\in\w}, (g_i)_{i\in\w})$
has the following properties:

\begin{itemize}
\setlength{\itemsep}{1mm}
\item[(a)]
$f_i \in G(N_i) \subset G(M_{2i+1})$ and
$g_i \in G(F_i) \subset G(M_{2i+1}) \subset G(M_{2i+2})$ for each $i \in \w$.

\item[(b)] $g_if_j = f_j g_i$ for $j \geq i+1$ since $F_i  \cap N_j = \emptyset$.
\item[(c)]
$h = \rho \sigma(h)
= \rho((f_i)_{i\in\w}, (g_i)_{i\in\w})
= \lambda_{\N}((f_i)_{i\in\w})\lambda_{\F}((g_i)_{i\in\w})
=  (f_0f_1f_2\cdots)(g_0g_1g_2\cdots)
= f_0g_0f_1g_1f_2g_2 \cdots.$
\item[(d)] $(\id_M, f_0, g_0, f_1, g_1, f_2, g_2, \dots) \in \cbox_{i\in\w} G(M_i)$ and
$h = p\big(\id_M, f_0, g_0, f_1, g_1, f_2, g_2, \dots\big)$.
\end{itemize}
\vskip 1mm
Finally the required local section at $\id_M$ of the map $p : \cbox_{i\in\w} G_i \to G$ is defined by
$$\gamma : \U \to\cbox_{i\in\w} G_i :  \ \gamma(h) = \big(\id_M, f_0, g_0, f_1, g_1, f_2, g_2, \dots\big).$$
This completes the proof.
\end{proof}

\begin{proof}[\bf Proof of Proposition~\ref{prop_H_i}]
(1)$'$ Since $H_i$ is a connected closed subgroup of $G$ and $H_i \subset H_{i+1}$ for each $i \in \w$,
the union $H' = \cup_{i \in \w} H_i$ is a connected subgroup of $G$.
Since $H_i$ is open in $G_i$ for each $i \in \w$, by Proposition~\ref{prop_G_i} (4)
the subgroup $H'$ is open (and closed) in $G$ and hence $H = H'$.

(2)$'$ The properties $(1) - (4)$ in Proposition~\ref{prop_G_i} are verified as follows:

(1)  Since $G$ is paracompact, so is the closed subspace $H$.
Since $G$ is not first countable, the open subspace $H$ is not first countable.

(2) (i) Since $G_i$ an $l_2$-manifold, so is the open subgroup $H_i$.

(ii) (iii) We apply Corollary~\ref{cor_bdl-2} (3)(ii) to $K_{i+1} \subset K_i$ to conclude that
$G_i \cap H_{i+1}$ is LTC in $H_{i+1}$ and
$H_{i+1}/(G_i \cap H_{i+1}) \approx \EE_{K_{i+1}}(K_i, M)_0$, which is an $l_2$-manifold by Theorem~\ref{thm_Frechet}.
Since $H_i$ is an open subgroup of $G_i \cap H_{i+1}$,
the quotient space $(G_i \cap H_{i+1})/H_i$ is a discrete space and hence $H_i$ is LTC in $G_i \cap H_{i+1}$.
Hence from Lemma~\ref{lem_ltc}\,(2)  it follows that $H_i$ is also LTC in $H_{i+1}$.
By Lemma~\ref{lem_ltc}\,(3) the projection $H_{i+1}/H_i \to H_{i+1}/(G_i \cap H_{i+1})$ is a locally trivial bundle with the discrete fiber $(G_i \cap H_{i+1})/H_i$. Since $H_{i+1}/(G_i \cap H_{i+1})$ is an $l_2$-manifold,
$H_{i+1}/H_i$ is locally homeomorphic to $l_2$.
Since $H_{i+1}/H_i$ is metrizable by Lemma~\ref{lem_ltc}\,(1), we conclude that $H_{i+1}/H_i$ is also an $l_2$-manifold.

(3) (4) Since each $H_i$ is an open subgroup of $G_i$,
the small box $\cbox_{i \in \w} H_i$ is an open neighborhood of $(\id_M)_{i \in \w}$ in $\cbox_{i \in \w} G_i$.
Thus any local section $\sigma$ of the map $p : \cbox_{i \in \w} G_i \to G$ at $\id_M$ with $\sigma(\id_M) = (\id_M)_{i \in \w}$ restricts to that of the left multiplication map $p : \cbox_{i \in \w} H_i \to H$.
The remaining assertions follow from  Remark~\ref{rmk_op}.
\end{proof}

For the proof of Corollary~\ref{cor_top_type} we need a preliminary.
For any pairs of spaces $(X, A)$ and $(Y, B)$ denote by $[X,A; Y, B]$
the set of homotopy classes $[g]$ of maps of pairs $g : (X, A) \to (Y, B)$.
Any map of pairs $f : (Y, B) \to (Z, C)$ induces a function
$$f_{\#} : [X, A; Y, B] \to [X, A; Z, C], \ f_{\#} : [g] \mapsto  [fg].$$

Suppose $L$ is a compact space and $K$ is a closed subset of $L$.
The inclusion maps $H_i \subset H_{i+1}$ and $H_i \subset H$ $(i \in \w)$
induce the associated functions between pointed sets:
$$\xymatrix@M+1pt{
[L, K; H_i, \id_M] \ar[rr]^{} \ar[dr]_{} & & [L, K; H_{i+1}, \id_M] \ar[dl]^{} \\
& [L, K; H, \id_M] &
}$$
Taking the direct limit, we obtain a function between pointed sets
$$\iota : \varinjlim\,[L, K ; H_i, \id_M] \longrightarrow [L, K; H, \id_M].$$
Since any compact subset of $H$ is included in some $H_i$, we have the following conclusion.

\begin{lemma} For any pair of compact spaces $(L, K)$ the inclusion induced function
$$\iota : \varinjlim\,[L, K ; H_i, \id_M] \longrightarrow [L, K; H, \id_M]$$
is a bijection.
\end{lemma}

For $m\in\w\cup\{\infty\}$, a map $f : X \to Y$ between path-connected spaces is called an {\em $m$-equivalence} if for some base point $x \in X$, the induced homomorphism on the $k$-th homotopy group
$$f_\# : \pi_k(X, x) \to \pi_k(Y, f(x))$$
is an isomorphism for all $k<m$ and an epimorphism for $k = m$.
An $\infty$-equivalence is called a {\em weak equivalence}.
If both $X$ and $Y$ have the homotopy type of CW-complexes, then every weak equivalence is a homotopy equivalence. Note that the groups $H$ and $H_i$ $(i \in \w)$ are path-connected and have the homotopy type of CW-complexes.

\begin{corollary}\label{cor_homotopy_type} Let $m\in\w\cup\{\infty\}$. If each inclusion $H_i \subset H_{i+1}$ is an $m$-equivalence, then so is the inclusion $H_1 \subset H$.
For example, if each $H_i$ is contractible, then so is $H$ and hence $H \approx l_2 \times \IR^\infty$.
\end{corollary}

\begin{proof}[\bf Proof of Corollary~\ref{cor_top_type}]
We keep the notations $M_i$, $K_i$, $H_i$, $i \in \w$, and $H$.

(1), (2) Since $M$ is connected, we may assume that
for each $i \in \w$ (a) $M_i$ is connected and
(b) each connected component of $K_i = M \setminus \Int M_i$ is non-compact.
By Corollary~\ref{cor_homotopy_type} it suffices to show that each $H_i$ is contractible.
Note that the restriction map $H_i \to \DD_0(M_i, \partial M_i) : h \mapsto h|_{M_i}$ is a homotopy equivalence.

For $n = 1 ,2$ the assertion follows from \cite{EE}, \cite[Section 2.7]{Ivanov3}, \cite{Sc}, \cite{Sm}, etc.
In the case $n = 3$, if $M_i$ is a 3-ball, then $\DD_0(M_i, \partial M_i)$ is contractible
 by the Smale conjecture \cite[Appendix (1)]{Hat2}.
If $M_i$ is not a 3-ball, then
 by the assumption, $M_i$ is an orientable Haken 3-manifold
 with boundary \cite{Hempel, Waldhausen} and $\DD_0(M_i, \partial M_i)$ is contractible by \cite{Hat1}, \cite{Ivanov1}, \cite{Ivanov2}.

(3) Take a collar $\partial N \times [0,1]$ of
$\partial N = \partial N  \times \{0\}$ in $N$ and
let $M_i = N \setminus (\partial N \times [0, 1/(i+1)))$ and
$K_i = M\setminus\Int M_i = \partial N \times (0, 1/(i+1)]$ for $i \in \w$.
First we show that the inclusion $H_i \subset H_{i+1}$ is a homotopy equivalence for each $i \in \w$.
The space $\EE_{K_{i+1}}(K_i, M)$ is contractible
since it is the space of embeddings of the collar $K_i$ relative to $K_{i+1}$.
In particular, $\EE_{K_{i+1}}(K_i, M) = \EE_{K_{i+1}}(K_i, M)_0$ and
by Corollary~\ref{cor_bdl-2}\,(2)(ii)
the restriction map $r : H_{i+1} \to \EE_{K_{i+1}}(K_i, M)$ is
a fiber bundle with fiber $H_{i+1}\cap G_i$.
Since the base space $\EE_{K_{i+1}}(K_i, M)$ is contractible and paracompact,
this bundle is trivial and there exists a fiber preserving homeomorphism
$$\phi : H_{i+1} \longrightarrow \EE_{K_{i+1}}(K_i, M) \times (H_{i+1}\cap G_i).$$
Since $H_{i+1}$ is connected, we have $H_i = H_{i+1}\cap G_i$.
Since $\phi$ is fiber preserving, we have the homeomorphism of pairs
$$\phi : (H_{i+1}, H_i) \longrightarrow (\EE_{K_{i+1}}(K_i, M), \{ i_{k_i} \}) \times H_i.$$
Since $\EE_{K_{i+1}}(K_i, M)$ is contractible,
the inclusion $\{ i_{k_i} \} \times H_i \subset \EE_{K_{i+1}}(K_i, M) \times H_i $ is a homotopy equivalence and
so is the inclusion $H_i \subset H_{i+1}$.

From Corollary~\ref{cor_homotopy_type} it follows that $H \simeq H_1 \simeq \DD_0(M_1; \partial M_1) \approx
\DD_0(N; \partial N)$. Since the last one is an $l_2$-manifold, we have
$H \approx \DD_0(N;\partial N) \times \IR^\infty$ by Corollary~\ref{cor_factor}.
\end{proof}

\end{document}